\documentclass[12pt]{article}
\usepackage{amsmath}
\usepackage{amsfonts}
\usepackage{amssymb}
\usepackage{latexsym}
\usepackage{euscript}
\usepackage{a4wide}

\providecommand{\bysame}{\leavevmode\hbox to3em{\hrulefill}\thinspace}
\def\lra{\longrightarrow}
\def\ra{\rightarrow}
\def\G{\Gamma}
\def\frak{\mathfrak}

\def\rank{\mathrm{rank}}

\def\text{\textrm}

\def\Mt{\M_f}
\def\Mf{\M_t}
\def\ord{\mathrm{ord}}
\def\GL{\mathrm{GL}}

\def\l{\ell}
\def\Qbar{\overline{\Q}}

\def\Gal{\mathrm{Gal}}
\def\F{\Bbb F}
\def\Rem{\noindent \bf Remark\rm. }

\def\Proof{\noindent \bf Proof\rm. }
\def\mod{\text{ mod }}

\def\qed{\hfill \square \ }
\def\qeda{\qed}

\def\Of{{\cal O}}
\def\C{\underline{C}}
\def\I{\mathcal I}
\def\D{\mathcal D}
\def\Di{\frak{D}}
\def\M{\mathcal M}
\def\ovM{\overline{\mathcal M}}
\def\ovMt{\ovM_f}
\def\ovMf{\ovM_t}
\def\haM{\widehat{\M}}
\def\haMt{\haM_f}
\def\haMf{\haM_t}
\def\haA{\hat{A}}

\def\Cl{\mathrm{Cl}}
\def\Bbb{\mathbb}
\def\Q{\Bbb Q}
\def\Z{\Bbb Z}
\def\T{\Bbb T}

\def\A{\mathcal A}
\def\l{\ell}
\def\GL{\mathrm{GL}}
\def\Aut{\mathrm{Aut}}
\def\mod{\text{ mod }}
\def\Ext{\mathrm{Ext}}
\def\Id{\mathrm{Id}}

\begin{document}
\newtheorem{sublemma}{Sub-lemma}
\newtheorem{theorem}{Theorem}[section]
\newtheorem{lemma}[theorem]{Lemma}
\newtheorem{df}[theorem]{Definition}
\newtheorem{cor}[theorem]{Corollary}

\author{Frank Calegari}
\title{Semistable abelian Varieties over $\Q$}
\maketitle
\abstract{We
prove that for $N=6$ and $N=10$, there do not exist any non-zero 
semistable abelian 
varieties over $\Q$ with good reduction outside primes dividing $N$.
Our results are contingent on the GRH discriminant bounds of
Odlyzko. Combined with recent results of Brumer--Kramer and of
Schoof, this result is best possible: if $N$ is squarefree,
there exists a 
non-zero semistable abelian variety over $\Q$ 
with good reduction outside primes dividing $N$
precisely when $N \notin \{1,2,3,5,6,7,10,13\}$.
\footnote{2000 Mathematics Subject Classification: $\mathbf{14K15}$}

\section{Introduction.}

In $1985$, Fontaine \cite{Fontaine} proved a conjecture
of Shafarevich to the effect that there do not
exist any
nonzero
abelian varieties over $\Z$
(or equivalently, abelian varieties $A/\Q$  with good reduction
everywhere). Fontaine's
approach was via finite group schemes over local fields.
In particular, he proved the following theorem:
\begin{theorem}[Fontaine] \label{theorem:font}  Let
$G_{\l}$ be a finite flat
group scheme over $\Z_{\l}$ killed by $\l$. Let
$L = \Q_{\l}(G_{\l}):=\Q_{\l}(G_{\l}(\Qbar_{\l}))$. Then
$$v(\frak{D}_{L/\Q_{\l}}) <  1 + \frac{1}{\l-1}$$
where $v$ is the valuation on $L$ such that $v(\l) = 1$,
and $\frak{D}_{L/\Q_{\l}}$ is the different of
$L/\Q_{\l}$. \end{theorem}
If $G_{\l}$ is the restriction of some finite flat
group scheme
$G/\Z$ 
 then $\Q(G)$ is
\emph{a fortiori} unramified at primes
outside $\l$. In this context, the result of Fontaine is striking since
it implies that the field $\Q(G)$ has particularly small
root discriminant.
If $A/\Q$ has good reduction everywhere, then it has a
smooth proper  N\'{e}ron model $\A/\Z$, and $G:=\A[\l]/\Z$ is a
finite flat group scheme.
Using the discriminant bounds of
Odlyzko \cite{Odlyzko},
Fontaine showed that for certain
small primes ${\l}$, for every $n$,
either $A/\Z$ or some isogenous abelian variety has
a rational ${\l}^n$-torsion point. Reducing $A$ modulo $p$
for some prime $p$ of good reduction (in this case, any prime),
one finds abelian varieties  over $\F_{p}$ with
at least ${\l}^n$ rational points. One knows, however, that isogenous
abelian
varieties over $\F_p$ have an equal and thus bounded
number of points. This contradiction proves that $A/\Q$ cannot exist.

If one considers
abelian varieties $A/\Q$ such that $A$ has
good reduction outside a
single prime $p$, one can no longer
expect nonexistence results. Indeed,
there exist abelian varieties with good reduction everywhere
except at $p$.
One such class of examples are the Jacobians of modular curves
$X_0(p^n)$, which have positive genus
for every $p$ and sufficiently large $n$.
A natural subclass of abelian varieties, however, are
the semistable ones. By considering the modular abelian varieties
$J_0(N)$, one finds nonzero
semistable abelian varieties over $\Q$ which have good
reduction outside $N$ for
all squarefree
$N \notin \{1,2,3,5,6,7,10,13\}$.
A reasonable conjecture to make is that there are no
semistable abelian varieties over $\Z[1/N]$ for
$N$ in this set.
Fontaine's Theorem is the case
$N = 1$.
Recently Brumer and
Kramer  \cite{BK}
 proved this statement for $N \in \{2,3,5,7\}$,
and (by quite different methods) Schoof \cite{Schoof}
 for $N \in \{2,3,5,7,13\}$.
In this paper, we treat the remaining cases $N \in \{6,10\}$,
and prove the following theorems:
\begin{theorem} \label{theorem:results}
Let $A/\Q$ be an abelian variety with everywhere semistable
reduction, and good reduction outside $2$ and $3$.
If
the GRH discriminant bounds of Odlyzko hold, then $A$ has
dimension $0$.
\end{theorem}

\begin{theorem} \label{theorem:results5}
Let $A/\Q$ be an abelian variety with everywhere semistable
reduction, and good reduction
outside $2$ and $5$. If 
the GRH discriminant bounds of Odlyzko hold, then $A$ has
dimension $0$.
\end{theorem}
We note that in Fontaine \cite{Fontaine},
Brumer--Kramer \cite{BK} and Schoof \cite{Schoof}, the
GRH is \emph{not} assumed.
Our technique for proving these results is  linked strongly
to the ideas in
Brumer--Kramer \cite{BK} and Schoof \cite{Schoof}, and
thus we consider it important to
briefly recall the main ideas of these papers now.
Schoof's
approach is similar in spirit to Fontaine's. Instead of working  with
finite flat group schemes over $\Z$,  one considers finite flat
group schemes
over $\Z[\frac{1}{p}]$, where $p$ is prime.
In order to restrict to the class of group schemes 
possibly arising from non-semistable
abelian varieties, one uses the following fact due to Grothendieck
(\cite{Groth}, Expos\'{e} IX, Prop.~3.5):

\begin{theorem}[Grothendieck] \label{theorem:Groth} Let $A$ be an
abelian variety over $\Q$ with semistable reduction at $p$. Let
$\I_p \subset \Gal(\Qbar/\Q)$ denote a choice of inertia group at $p$. Then the
 action of $\I_p$ on the $\l^n$-division
points of $A$ for $\l \ne p$
is rank two unipotent; i.e., as an endomorphism, for
$\sigma \in \I_p$,
$$(\sigma - 1)^2 A[\l^n] = 0.$$
In particular, $\I_p$ acts through its maximal pro-$\l$ quotient,
which is procyclic.
\end{theorem}
Thus one may
 restrict attention to finite flat
group schemes
$G/\Z[\frac{1}{p}]$  of $\l$-power order such that inertia at $p$
acts through its maximal pro-$\l$ quotient.
The key step of Schoof's
 approach is to show that any such group scheme
admits a filtration by the group schemes $\Z/\l \Z$ and $\mu_{\l}$.
Using this filtration, along with various extension results (in the
spirit of Mazur \cite{Mazur}, in particular
Proposition 2.1 pg. 49 and Proposition 4.1 pg. 58) for
group schemes over $\Z[\frac{1}{p}]$, one shows as in Fontaine that
for each $n$,
some abelian variety isogenous to $A$ has rational
torsion points of order $\l^n$.
The approach of Brumer and
Kramer is quite different. Although, as in Schoof and Fontaine, they
use discriminant bounds to control $\Q(A[\l])$ for particular $\l$, they
seek a contradiction not to any local bounds but
to a theorem of Faltings.
Namely, they construct
infinitely many pairwise non isomorphic but isogenous 
abelian varieties, contradicting
the finiteness of this set (as follows from Faltings
\cite{Faltings}, Satz 6,
pg. 363).
The essential difference in the two
approaches, however, is that 
Brumer and Kramer use the explicit  description of
the Tate module $\T_{\l}$ of $A$ at a prime $p$ of semistable reduction.
Such a description is once more due to Grothendieck \cite{Groth}.

Both of these approaches fail (at least na\"{\i}vely) to work 
when $N = 6$ or $10$. Using Schoof's approach
one runs into a problem (when $N = 10$)
 because $\mu_3$ admits  non-isomorphic
finite flat group scheme
extensions by $\Z/3 \Z$ over $\Z[\frac{1}{10}]$, whereas no nontrivial
 extensions
exist over either $\Z[\frac{1}{2}]$ or $\Z[\frac{1}{5}]$. One
difficulty that arises in
Brumer and Kramer's approach is that
the field $\Q(A[3])$  fails to
have a unique prime above the bad primes $2$ or $5$,  as fortuitously
happens in the cases they consider.
We combine both methods, as well as some new ideas, to prove
our results.
In the next section we recall some
definitions and results from Brumer and Kramer's
paper.

\subsection{Notation.}

Let $p$ be a 
prime number. Let $\D_p = \Gal(\Qbar_p/\Q_p)$ denote
the local Galois group at~$p$. For a Galois
 extension of global fields $L/\Q$, we
denote a decomposition group at $p$ by $\D_p(L/\Q)$. This is well
defined up to conjugation, or equivalently, up to an embedding 
$L \hookrightarrow \Qbar_p$ which we shall fix when necessary. 
In the same spirit, let $\I_p
= \Gal(\Qbar_p/\Q^{unr}_p)$, and let $\I_p(L/\Q)$ be
an inertia group at $p$ as a subgroup of $\D_p(L/\Q)$ and of
$\Gal(L/\Q)$. One notes that $\I_p$
is normal in $\D_p$. For any $\D_p$-module $\M$,
let $\ovM$ denote $\M/\ell \M$; it is a $\D_p$-module
killed by $\ell$. We shall use  $\haM$ to
denote a $\Gal(\Qbar/\Q)$-module killed by $\ell$
constructed functorially
from $\ovM$.
A ``finite'' group scheme $G/R$ will always mean a group scheme
$G$ finite and flat over $\mathrm{Spec} \ R$.
For an abelian variety $A$, let $\haA$ denote
the dual abelian variety.

\section{Local Considerations.}

\subsection{Preliminaries.}

In this section we  introduce some notation and results from
the paper of Brumer and Kramer \cite{BK}.

Let $A/\Q$ be an abelian variety of dimension $d > 0$ with
semistable reduction at $p$.
Let $\l$ be a prime different from $p$,
and consider the Tate module
$\T_{\l}(A/\Q_p)$. 
Let $\Mt(p) = \T_{\l}(A/\Q_p)^{\I_p}$,
and let $\Mf(p)$ be the submodule of $\T_{\l}(A/\Q_p)$ orthogonal
to  $\Mt(p)(\hat{A})$ under the Weil paring
$$e: \T_{\l}(A) \times \T_{\l}(\hat{A}) \lra \Z_{\l}(1).$$
In Brumer and Kramer, these modules were referred to
as $\M_1$ and $\M_2$ respectively. Following a suggestion,
we use instead the hopefully more suggestive notation
$\Mt$ ($f$ for finite or fixed) and $\Mf$ ($t$ for toric).
Since $A$ is semistable, there are inclusions
$$0 \subseteq \Mf(p) \subseteq \Mt(p) \subseteq \T_{\l}(A/\Q_p).$$
Since $\I_p$ is normal in $\D_p$, the groups
$\Mt(p)$ and $\Mf(p)$ are $\D_p = \Gal(\Qbar_p/\Q_p)$-modules.
Let $\A/\Z$ be the N\'{e}ron model for $A$. Let $\A^0_{\F_p}$
be the connected component of the special fibre of $\A$ at $p$.
It is an extension of an abelian variety of dimension $a_p$ by
a torus of dimension $t_p = d - a_p$. One has
$\rank(\Mf(p)) = t_p$ and $\rank(\Mt(p)) = t_p + 2 a_p = d + a_p$.

\begin{df}[Brumer--Kramer] Let $A$ be an abelian
variety with bad reduction at $p$. 
Let $i(A,\l,p)$ denote the 
minimal integer $n \ge 1$ such that
$\Q_p(A[\l^n])$ is ramified at $p$.
Call $i(A,\l,p)$ the ``effective stage of inertia''. 
\end{df}

We note that
$i(A,\l,p)$  is finite
by the criterion of N\'{e}ron--Ogg--Shafarevich.

Let $\Phi_{A}(p) = (\A/\A^0)(\overline{\F}_p)$ be the component
group of $A$ at $p$. For a finite group $G$, let $\ord_{\l}(G)$
denote the largest exponent $d$ such that $\l^d$ divides the
order of $G$.
Recall the following result from \cite{BK}:

\begin{theorem}[Brumer--Kramer] \label{theorem:BK} 
Let $A$ be a semistable abelian variety with bad reduction at $p$.
 Let
$\ovMt(p)$ and $\ovMf(p)$ denote the projections
of $\Mt(p)$ and $\Mf(p)$ to $A[\l]$. 
Suppose that $\kappa$ is a  $\Gal(\Qbar_p/\Q_p)$-submodule
of $A[\l]$ and let $\phi:A \lra A'$ be a $\Q_p$-isogeny
with kernel $\kappa$. Then
$$\ord_{\l}(\Phi_{\haA'}(p)) - \ord_{\l}(\Phi_{\haA}(p))
= \dim ( \kappa \cap \ovMf(p)) + \dim ( \kappa \cap \ovMt(p)) - 
\dim \ \kappa.$$
Moreover, if $\ovMf(p) \subseteq \kappa \subseteq \ovMt(p)$, then
$i(A',\l,p) = i(A,\l,p) + 1$.
\end{theorem}

By taking $\kappa$ to be a proper $\Gal(\Qbar/\Q)$ submodule
of $A[\l]$,
Brumer and Kramer use this theorem  to construct infinitely
many non-isomorphic varieties isogenous to $A$ over $\Q$.
This contradicts
Faltings' Theorem.  Although we shall also use Faltings' Theorem,
our final contradiction will come from showing that $A$ (or some isogenous
abelian
variety) has too many points over some finite field, contradicting 
Weil's Riemann hypothesis, much as in the approach of
Schoof \cite{Schoof}.
We shall also make use of the following lemma.
\begin{lemma} \label{lemma:sigmaminus} Let $\sigma \in \I_p$.
The image of $(\sigma - 1)$ acting
on $\T_{\l}(A)$ lies in $\Mf(A)$. The image of
$(\sigma-1)$ on $A[\l]$ lies in $\ovMf(p)$.
\end{lemma} 
\Proof Let $y \in \Mt(p)(\haA)$, and $x \in \T_{\l}(A)$.
Then $e((\sigma -1)x,y) = e(x^{\sigma},y)/e(x,y)$. Since
both $y$ and $\Z_{\l}(1)$ are fixed by $\sigma$, we conclude
that 
$$e((\sigma -1)x,y) = e(x^{\sigma},y^{\sigma})/e(x,y)
= e(x,y)^{\sigma}/e(x,y) = 1.$$
Thus $({\sigma-1})x \in \Mf(p)$.  The second statement
of Lemma~\ref{lemma:sigmaminus} follows from the first. $\qeda$

\subsection{Results.}
\label{sec:results}

In proving Theorem~\ref{theorem:results} (or \ref{theorem:results5}),
we may assume that $A$ has bad reduction at both $2$ and $3$ 
(respectively, both $2$ and $5$), since otherwise we may apply
the previous
results of Brumer--Kramer \cite{BK}, Schoof \cite{Schoof},
or Fontaine \cite{Fontaine}.

The proof of Theorem~\ref{theorem:results5} is very similar
to the proof of Theorem~\ref{theorem:results},
 although some
additional complications arise. Thus we restrict ourselves
first to the case $N = 6$, and then later explain how our
proof  can be adapted to work for $N = 10$. 
One main ingredient is the following result, proved in
section~\ref{sec:groups}: 

\begin{theorem} \label{theorem:5power} 
Let $G/\Z[\frac{1}{6}]$ be a finite group scheme of \ $5$-power
 order such that:
\begin{enumerate}
\item  Inertia at $2$ and $3$ acts through
a cyclic $5$-group.
\item  The action of inertia on the subquotients
$G[5^n](\Qbar)/G[5^{n-1}](\Qbar)$ is through an order $5$ quotient for all $n$.
\end{enumerate}
Assume the GRH
discriminant bounds of Odlyzko. Then $G$ has a filtration by
the group schemes $\Z/5 \Z$ and $\mu_5$. Moreover, if $G$ is
killed by $5$, then $\Q(G) \subseteq K$,
where $K:= \Q(\sqrt[5]{2},\sqrt[5]{3},\zeta_5)$.
\end{theorem}

In particular,
if $A/\Q$ is a semistable
abelian variety with good reduction outside $2$ and $3$, and
$\A/\Z$ is its  N\'{e}ron model, then by Theorem~\ref{theorem:Groth}
the conditions of Theorem~\ref{theorem:5power} are satisfied by
the finite group scheme  $\A[5^n]/\Z[\frac{1}{6}]$ for each $n$.
Thus $\A[5^n]$ has a filtration by
the group schemes $\Z/5 \Z$ and $\mu_5$, and $\Q(A[5]) \subseteq K$.
These results and their proofs are
 of  the same flavour as results in Schoof \cite{Schoof}.
One such result from that paper
 we  use explicitly is the following (a special case of Theorem
3.3 and the proof of Corollary 3.4 in {\it loc. cit.}\rm):

\begin{theorem}[Schoof]  \label{theorem:Schoof} 
Let $p = 2$ or $3$. 
Let $G/\Z[\frac{1}{p}]$ be a finite group scheme of
$5$-power order such that inertia at $p$ acts through
a cyclic $5$-group. Then $G$ has a filtration by
the group schemes $\Z/5 \Z$ and $\mu_5$. Moreover, the extension
group $\Ext^1(\mu_5,\Z/5 \Z)$ of group schemes over $\Z[\frac{1}{p}]$ is
trivial, and there exists an exact sequence of group schemes
$$0 \lra M \lra G \lra C \lra 0$$
where $M$ is a diagonalizable group scheme over $\Z[\frac{1}{p}]$,
 and $C$ is
a constant group scheme.
\end{theorem}

In sections~\ref{sec:const}, \ref{sec:nontoric} and \ref{sec:toric} 
we shall assume 
there exists a semistable abelian variety  $A/\Z[\frac{1}{6}]$,
and derive a contradiction using Theorem~\ref{theorem:5power}.

\subsection{Construction of Galois Submodules.}
\label{sec:const}

The proof of Brumer and Kramer relies on the fact that for
abelian varieties with bad semistable reduction at one prime
$p \in \{2,3,5,7\}$, there exists an $\l$ such that
there is a unique prime above $p$ in $\Q(A[\l])$.
In this case, the $\D_p$-modules $\ovMt(p)$ and
$\ovMf(p)$ are  automatically $\Gal(\Qbar/\Q)$-modules, and so
one has a source of $\Gal(\Qbar/\Q)$-modules with which
to apply Theorem~\ref{theorem:BK}.
This approach fails in our case, (at least if $\l = 5$) since
Theorem~\ref{theorem:5power}  allows the possibility that
$\Q(A[5])$ could be as big as
$K:= \Q(2^{1/5},3^{1/5},\zeta_5)$, and   $2$ and $3$ split into
$5$ distinct primes in $\Of_K$.
On the other hand, something
fortuitous does happen, and that is that the inertia subgroups
$\I_p(K/\Q)$ for $p = 2$, $3$ are \emph{normal} subgroups of
$\Gal(K/\Q)$, when \emph{a priori} they are only normal
subgroups of $\D_p(K/\Q)$. Using this fact we may 
construct global Galois modules from the local $\D_p(K/\Q)$-modules
$\ovMt(p)$ as follows.

\begin{lemma} \label{lemma:mod}   
Let $F = \Q(A[\l])$,
$\G = \Gal(F/\Q)$,
and  $H \subseteq \G$ be a normal subgroup of $\G$.
Let $\ovM$ be a subgroup of $A[\l]$ fixed pointwise by $H$.
Let $\haM$ be the
$\Gal(\Qbar/\Q)$-submodule generated by the points of $\ovM$. Then
$\Q(\haM) \subseteq E$, where $E$ is the
fixed field of $H$.
\end{lemma}

 \Proof
By Galois theory, it suffices to show that $\haM$ is fixed by $H$.
This result is a special case of the more general fact: 
If $H$ is any normal subgroup of $\G$, then
any $\G$-module generated by $H$-invariant elements  is itself
$H$-invariant.
Any sum  of elements fixed by $H$ is clearly fixed
by $H$. Thus it remains to show that any element  $gP$
with $g \in \G$ is also fixed by $H$. For
this we observe that
$$h(gP) = g(g^{-1}hgP) = gP$$
since $g^{-1}hg \in H$. $\qed$

\

\begin{df} Let
$\haMt(p)$  be the $\Gal(\Qbar/\Q)$-module generated  by
$\ovMt(p)$, considered as a
subgroup of $A[\l]$ after choosing some embedding $\Qbar
\hookrightarrow \Qbar_p$. 
\end{df} 

Since all embeddings $\Qbar \hookrightarrow \Qbar_p$
differ by an automorphism of $\Qbar$, we find that $\haMt(p)$ does
not depend on the choice of embedding, although $\ovMt(p)$
does, in general. 
We note that by Faltings theorem, there exist only finitely
many abelian varieties over $\Q$ isogenous to $A$. Thus is
makes sense to chose a representative from the isogeny class of $A$
that is \emph{maximal} with respect to any well defined property.

\begin{lemma} \label{lemma:newreview}
Suppose that
$\ord_{5}(\Phi_{\haA}(2))$ is maximal amongst all abelian varieties isogenous to $A$. 
Then
\begin{enumerate}
\item $A[5]$ is unramified at $2$
\item There is an exact sequence
$$0 \lra \mu^m_5 \lra A[5] \lra (\Z/5 \Z)^n \lra 0$$
with $m + n = 2d$. Moreover, $m=n=d$.
\end{enumerate}
Similarly, if $A$ is chosen such that
$\ord_{5}(\Phi_{\haA}(3))$  is maximal, then $A[5]$ is unramified
at $3$ and statement $2$ still holds. Finally, $A$ and any variety isogenous
to $A$ has ordinary reduction at $5$. 
\end{lemma}

%% A is ordinary at 5 regardless
%%
%%

Since $\I_p(K/\Q)$ is a normal subgroup of $\Gal(K/\Q)$,
Lemma~\ref{lemma:mod} implies that
$$\Q(\haMt(2)) \subseteq \Q(\zeta_5,3^{1/5}), \qquad
\Q(\haMt(3)) \subseteq \Q(\zeta_5,2^{1/5}).$$
We now apply Theorem~\ref{theorem:BK} with
$\kappa = \haMt(2)$. Let $A' = A/\kappa$. 
Since $\kappa$ is a $\Gal(\Qbar/\Q)$ module $A'$
is an abelian variety over $\Q$. We see that
$$\ord_{5}(\Phi_{\haA'}(2)) - \ord_{5}(\Phi_{\haA}(2))
=  \dim \ \kappa \cap \ovMf(2)
+ \dim \ \kappa \cap \ovMt(2) - \dim \ \kappa.$$
Since by construction $\ovMf(2) \subseteq \ovMt(2)
 \subseteq \kappa$, the right hand side is equal to
$$2d - \dim \ \kappa \ge 0.$$
Yet from the maximality of $\ord_{5}(\Phi_{\haA}(2))$, it follows that
$2d - \dim \ \kappa \le 0$. Thus $\dim \ \kappa = d$, and in 
particular $\haMt(2) = \kappa = A[5]$.  Thus by Lemma~\ref{lemma:mod}  $A[5]$ is
unramified at $2$.
Note that this same construction can be applied
\emph{mutatis mutandis} when $2$ is replaced by $3$.
Since $A[5]$ is unramified at $2$,
it follows from a standard patching argument
(\cite{Mazur}, 1.2(b), p. 44) that
$A[5]$ prolongs to  a finite group scheme over $\Z[\frac{1}{3}]$.
Thus we may now apply  Theorem~\ref{theorem:Schoof}, and
conclude that there exists an exact sequence of
group schemes over $\Z[\frac{1}{3}]$
$$0 \lra \mu^m_5 \lra A[5] \lra (\Z/5 \Z)^n \lra 0$$
where $m + n = 2d$. 
It now remains to show that $m=n=d$.

Let $A' = A/\mu^m_5$.
The morphism $A \ra A'$ induces a proper map
$(\Z/5\Z)^n = A[5]/\mu^m_5 \ra A'$. By an fppf abelian sheaf argument,
we see that this map is a categorical monomorphism and
hence by EGA $\mathrm{IV}_3$ 8.11.5 (\cite{EGA}) a closed immersion. 
Specializing to the
fibre over $\F_5$ we find  that 
$$(\Z/5 \Z)^n \hookrightarrow A'_{\F_5}[5].$$
The $p$-rank of the  $p$-torsion subgroup  of an abelian variety 
over an algebraically closed field of
characteristic $p$ is at most the 
dimension $d$, with equality  only if 
$A$ is ordinary at $p$. Thus $n \le d$. Applying the same
argument to $\haA$ we find that $m \le d$ and thus
$n = m = d$, and $A$ has ordinary reduction at $5$. 
Since ordinary reduction is preserved under isogeny, we are done. $\qed$

\

We now divide our proof by contradiction into two cases.
In the first case we assume that
$A$ has mixed reduction at at least one of $2$ or $3$
(i.e. the connected component of the special fibre is
the extension of a \emph{non-trivial} abelian variety of
dimension
$a_p \ne 0$ by a torus of dimension $t_p=d-a_p$).
In the second case
we assume that $A$ has purely toric reduction at both $2$ and $3$.

\subsection{$A$ has Mixed Reduction at $2$ or $3$.}
\label{sec:nontoric}

Let  $\ord_{5}(\Phi_{\haA}(2))$ be maximal. Then from
Lemma~\ref{lemma:newreview} 
there is an exact sequence over $\Z[\frac{1}{3}]$
$$0 \lra \mu^d_5 \lra A[5] \lra (\Z/5 \Z)^d \lra 0.$$
If $A$ has mixed reduction at $2$ then
$a_2 > 0$, and  $\Mt(2)$ has
rank $t_2 + 2 a_2 = d + a_2 > d$. In particular,
$\kappa:= \ovMt(2) \cap \mu^d_5$ is nontrivial and
defines a  diagonalizable
$\Gal(\Qbar/\Q)$-submodule of $A[5]$
(here we use the fact that every subgroup of
$\mu^d_5(\Qbar)$ is $\Gal(\Qbar/\Q)$ stable). 
We now apply Theorem~\ref{theorem:BK}. Let $A' = A/\kappa$. We find that
$$\ord_{5}(\Phi_{\haA'}(2)) - \ord_{5}(\Phi_{\haA}(2))
=  \dim \ \kappa \cap \ovMf(2)
+ \dim \ \kappa \cap \ovMt(2) - \dim \ \kappa.$$
Since $\kappa \subseteq \ovMt(2)$, the last two terms cancel,
and  $\ord_{5}(\Phi_{\haA'}(2))$ is also maximal. Hence we may
repeat this process, thereby constructing morphisms
$A \lra A^{(n)}$ with larger and larger kernels $\kappa_n$, where
$\kappa_n$ has a filtration by copies of the finite group scheme
$\mu_5$. 

\begin{lemma} \label{lemma:mult}
Any extension of diagonalizable group schemes of $5$-power order
over $\Z[\frac{1}{6}]$ is diagonalizable.
\end{lemma}

\Proof By taking Cartier duals, it suffices to prove the
analogous statement for constant group schemes: Any extension of
$5$-power order constant group schemes over $\Z[\frac{1}{6}]$
is constant. The action of $\Gal(\Qbar/\Q)$ on any such extension
is unramified outside $2$ and $3$, and
acts via a $5$-group. Since $p$-groups are solvable,
it suffices to prove that there are no Galois $5$-extensions
of $\Q$ unramified outside $2$ and $3$.
Easy class field theory (for example, the Kronecker--Weber
theorem) shows that no such extensions exist.
 $\qed$

\

Thus we have proven that
for all $n$, there exist exact sequences
$$0 \lra \kappa_n \lra A[5^{k(n)}] \lra M_n \lra 0.$$
where $\kappa_n$ is a diagonalizable group scheme,
$k(n)$ the smallest integer such that $5^{k(n)}$ kills
$\kappa_n$,
and $M_n$ is the cokernel. 
Hence the variety
$\hat{A}/M^{\vee}_n$ contains the arbitrarily large constant group scheme
$\kappa^{\vee}_n$, and so, after choosing some
auxiliary prime $q$ of good reduction, we see that
 $(\hat{A}/M^{\vee}_n)(\F_q)$ can be 
arbitrarily large.
This contradicts the uniform boundedness of the number
of points over $\F_q$ for all varieties isogenous to $\hat{A}$ 
(indeed, the number of points for all such varieties
is equal).

If $A$ does not have purely toric reduction at $3$, a similar
argument applies.

\subsection{$A$ has Purely Toric Reduction at $2$ and $3$.}
\label{sec:toric}

Under this assumption, for $p \in \{2,3\}$, we have
 $\Mf(p) = \Mt(p)$, and so
we write both as $\M(p)$. Again
we assume that $\ord_{5}(\Phi_{\haA}(2))$ is maximal. In particular,
we may assume that $\haM(2) = A[5]$, that $\Q(A[5])$ is
contained in $\Q(\zeta_5,3^{1/5})$, and that we have an
exact sequence of group schemes over $\Z[\frac{1}{3}]$:
$$0 \lra \mu^d_5 \lra A[5] \lra (\Z/5 \Z)^d \lra 0.$$
By abuse of notation we may also think of this as an
exact sequence of $\Gal(\Qbar/\Q)$-modules.

\begin{lemma} \label{lemma:hat3}  The Galois modules
$\haM(3)$ and $\mu^d_5$ coincide. Equivalently, there is an equality
of Galois modules: $\haM(3) = \mu^d_5$.
\end{lemma}

\Proof 
First we show that $\ovM(2) \cap \mu^d_5 = \{0\}$.
If not, then since $\dim \ \ovM(2) = d$,  the module
$\ovM(2)$ would not surject onto $(\Z/5 \Z)^d$, 
and the elements of $\ovM(2)$ could
not possibly generate $A[5]$ as 
a $\Gal(K/\Q)$-module\footnote{Another  way to reduce
to the case where  $\ovM(2) \cap \mu^d_5 = \{0\}$ is as follows:
if this intersection was nontrivial,
we could  take
quotients repeatedly until the resulting
intersection \emph{was} trivial. If this process repeated
indefinitely,
 we could apply the  arguments of 
 section~\ref{sec:nontoric} to produce a contradiction.}.
Thus by dimension considerations, as a $\F_5$-vector space,
we have that
$A[5] = \mu^d_5 \oplus \ovM(2)$.

\

Let $L:=\Q(\zeta_5,3^{1/5})$. Then as we have noted,
$\Q(A[5]) \subseteq L$.
Let $\sigma$ generate 
$\I_3(L/\Q)$.
From Grothendieck's Theorem (Theorem~\ref{theorem:Groth}), 
we have $(\sigma-1)^2 = 0$
as an endomorphism
on $A[5]$. Thus $\ovM(2)+\sigma \ovM(2)$ is a well defined
$\I_3$-module. On the other hand, $\ovM(2)$ is a 
$\D_2(L/\Q)$-module, and  
$\I_3$ is a set of representatives for the left cosets of
$\D_2(L/\Q)$ in $\Gal(L/\Q)$.
Thus
$\ovM(2)+\sigma \ovM(2)$ is a $\Gal(L/\Q)$-module, and so
$$\haM(2) =  \ovM(2)+\sigma \ovM(2).$$
 Since $\dim_{\F_5} \haM(2) = \dim_{\F_5}(A[5]) = 2d$, 
by dimension considerations one must have $\sigma \ovM(2) \cap \ovM(2) = 0$.

\

The decomposition group of $L:=\Q(\zeta_5,3^{1/5})$ at $3$ is
the entire Galois group $\Gal(L/\Q)$, 
and the inertia group
$\I_3$ is equal to $\langle \sigma \rangle$.
We  show
that $\mu^d_5 \subset \haM(3)$ and $\haM(3) \subset \mu^d_5$.

Since $\sigma \ovM(2) \cap \ovM(2) = \{0\}$, we have
$\mathrm{ker}(\sigma - 1) \cap \ovM(2) = \{0\}$. An element
killed by $\sigma -1$ is exactly fixed by $\I_3$. Thus
the only elements of $A[5]$ fixed by $\I_3$ are those
in $\mu_5^d$.
Since (by Lemma~\ref{lemma:mod}) 
$\haM(3)$ is unramified at $3$,  we have
$\haM(3) \subseteq \mu^d_5$. On the other hand, $|\haM(3)|
\ge|\ovM(3)| = 5^d = |\mu^d_5|$.  Thus we are done. $\qed$

\

We now apply Theorem~\ref{theorem:BK} again with
$\kappa = \haM(3) = \mu^d_5$.
If $A' = A/\mu^d_5$, then since $\ovM(3) = \haM(3)$,
we have $i(A',5,3) = i(A,5,3) + 1 \ge 2$. On the other hand, 
we see from the exact sequence for $A[5]$ 
that $(\Z/5 \Z)^d \subset A'[5]$. 
By Theorem~\ref{theorem:5power}
and the proof of Lemma~\ref{lemma:newreview}  
 we infer that
there exists an exact sequence of group schemes
over $\Z[\frac{1}{6}]$:
$$0 \lra (\Z/5 \Z)^d \lra A'[5] \lra \mu^d_5 \lra 0.$$
Replace $A$ by $A'$. Since $\Q(A[5])$ is unramified
at $3$, we know that is must be contained within 
$\Q(\zeta_5,2^{1/5})$. Since $A$ is ordinary at $5$,
however, we may prove more. 

\begin{lemma} \label{lemma:above} The field $\Q(A[5])$ is
$\Q(\zeta_5)$.
There is only one prime above $3$ in the extension
$\Q(A[5])/\Q$.
\end{lemma}

\Proof Consider the action of $\I_5$ on $A[5]$. 
By Lemma~\ref{lemma:newreview}, $A$ is ordinary at $5$.
Thus  $A[5]$ as an $\I_5$-module is an extension of a
constant module of rank $d$ by
a cyclotomic module of rank $d$. The $(\Z/5 \Z)^d$ inside $A[5]$
 must intersect
trivially with this cyclotomic module. Thus it provides a splitting
of $A[5]$ as an $\I_5$-module into a product of a cyclotomic module 
and a constant module. Thus $\Q_5(A[5])$ is unramified over $\Q_5(\zeta_5)$.
The maximal extension of $\Q(\zeta_5)$ inside $K$ unramified at
$1 - \zeta_5$ is $\Q(\zeta_5,18^{1/5})$. 
Since $\Q(A[5])$ is also unramified over $3$
(as $i(A,5,3) \ge 2$),
$\Q(A[5])$ must be
exactly $\Q(\zeta_5)$. The second statement of
the lemma  clearly follows
from the first. $\qed$ 

\

The second part of Lemma~\ref{lemma:above} implies that
$\ovM(3)$ is a $\Gal(\Qbar/\Q)$-module, 
as in \cite{BK}. Applying
Theorem~\ref{theorem:BK} once more, with
 $\kappa = \haM(3) = \ovM(3)$, and setting
$A'=A/\kappa$, we find that
$$i(A',5,3) = i(A,5,3) + 1 \ge 3.$$
Replace $A$ by $A'$.
In particular, $\Q(A[5^2])$ is unramified at $3$. Thus
by Theorem~\ref{theorem:Schoof} there exists a filtration
of group schemes over $\Z[\frac{1}{2}]$:
$$0 \lra M \lra A[5^2] \lra C \lra 0$$
where $M$ is a diagonalizable group scheme, and $C$ is a constant 
group scheme.
Let $q \in \Z$ be a prime of good reduction. We observe that
the varieties $A/M$ and $\hat{A}/C^{\vee}$ contain constant
subgroup schemes of order $\# C$ and $\# M$ respectively.
It follows from Weil's Riemann Hypothesis
that abelian varieties of dimension $d$ over $\F_q$ have
at most $(1 + \sqrt{q})^{2d}$ points.
Thus $\#C \le (1 + \sqrt{q})^{2d}$ and $\#M \le (1 + \sqrt{q})^{2d}$, and
$$5^{4d} = \#A[5^2] = \#C \#M \le (1 + \sqrt{q})^{4d}.$$
Choosing $q = 7$, say, then since  $5 > 1 + \sqrt{7}$, we have
a contradiction if $d > 0$.
This completes the proof of Theorem~\ref{theorem:results} except
for Theorem~\ref{theorem:5power}, which we prove now.

\section{Group Schemes over $\Z[1/6]$.}
\label{sec:groups}

First, some preliminary remarks on group schemes. Here
we  follow  Schoof \cite{Schoof}.

\

Let $(\l,N) = 1$. Let $\C$ be the category of finite group
schemes $G$ over $\Z[1/N]$ satisfying the following properties:
\begin{enumerate}
\item $G$ is killed by $\l$: $G = G[\l]$.
\item  For all $p | N$, 
the action of $\sigma \in \I_p$ on $G(\Qbar_{p})$ is
either trivial or cyclic of order $\l$.
\end{enumerate}
For example, $\Z/\l \Z$ and $\mu_{\l}$ are objects of $\C$.
As remarked in \cite{Schoof}, this category is closed under
direct products, flat subgroups and flat quotients. Thus, to prove that
any object of $\C$ has a filtration by $\Z/\l\Z$ and
$\mu_{\l}$ it suffices to show that the only simple objects of
$\C$ are
$\Z/\l \Z$ and $\mu_{\l}$. If $A/\Q$ is a semistable abelian
variety with good reduction at primes not dividing $N$,
 then from Theorem~\ref{theorem:Groth}, we have $\A[\l] \in \C$.
Another class of examples are the group
schemes $G_{\epsilon}$ defined by  Katz--Mazur
(\cite{Katz} Chapter 8, Interlude 8.7, \cite{Schoof}); for any
unit $\epsilon \in \Z[1/N]$ they construct a group scheme
$G_{\epsilon} \in \C$ of order $\l^2$ killed by $\l$.
$G_{\epsilon}$ is an extension of
$\Z/\l \Z$ by $\mu_{\l}$, and
$G_{\epsilon}(\Qbar) =
\Q(\zeta_{\l},\epsilon^{1/\l})$. 

\

Let $N = 6$ and $\l = 5$.
To prove that the only simple objects of $\C$ are
$\mu_{5}$  and $\Z/5 \Z$, it  suffices to show that
the $\Qbar$ points of any
object of $\C$ are defined over the field $K$,  where
$K = \Q(\zeta_5,2^{1/5},3^{1/5})$, because of the
following result:

\begin{lemma} \label{lemma:simple}
Let $G/\Z[1/N]$ be a simple group scheme killed by $\l$,
where $(N,\l) = 1$. Let $L = \Q(G(\Qbar))$ and 
suppose that 
$\Gal(L/\Q(\zeta_{\l}))$ is an $\l$-group. Then $G$ is either
$\Z/\l \Z$ or $\mu_{\l}$.
\end{lemma}

\Proof Since any $\l$-group acting on $(\Z/\l \Z)^d$ has at least one 
(in fact $\l-1$) nontrivial fixed point, 
there exists a point $P$ of $G$ defined over
$\Q(\zeta_{\l})$. Since $G$ is simple, $P$ generates $G$ as
a Galois module and thus $\Q(G) \subseteq \Q(\zeta_{\l})$.
In particular, $G$ is unramified outside $\l$, hence 
unramified at each prime dividing $N$. Thus 
$G$ prolongs to a finite group scheme over 
$\Z$, killed by $\l$, and with $\Q(G) \subseteq \Q(\zeta_{\l})$. Since
the $(\l-1)$th roots of unity are in $\F^*_{\l}$, 
 and since all irreducible
$\F_{\l}$ representations of the abelian group
$\Gal(\Q(\zeta_{\l})/\Q)$ are one dimensional, 
any simple subgroup scheme of $G$ has order $\l$.
From Oort--Tate \cite{Oort}, the finite group schemes of
order $\l$ over $\Z$ are $\Z/\l \Z$ and $\mu_{\l}$. $\qed$

\

Let $G$ be an object of $\C$. To prove that
$\Q(G) \subseteq K$ it clearly suffices to prove the same inclusion 
for any group scheme which contains $G$ as a direct
factor. 
Consider the field
$L = \Q(G \times G_{-1} \times G_{2} \times G_{3})$.
One sees from the definition of $G_{\epsilon}$
that $K:=\Q(\zeta_5,2^{1/5},3^{1/5}) \subseteq L$.
We prove that $L = K$.
Using the estimates  of
Fontaine \cite{Fontaine} we obtain an upper bound on the ramification
of $L$ at $5$. Since inertia at $2$ and $3$ acts through a
cyclic subgroup of order $5$, we also have ramification
bounds at $2$ and $3$. As in 
Schoof \cite{Schoof} and Brumer--Kramer \cite{BK},
we obtain   the following estimate
of the root discriminant
$$\delta_L:=|\Delta_L|^{1/[L:\Q]} <  5^{1 + \frac{1}{5-1}} 2^{1 - \frac{1}{5}}
3^{1 - \frac{1}{5}} = 5^{5/4} 6^{4/5} = 31.349 < 31.645.$$
Note that this is not only an inequality of real numbers:
one also has that as algebraic integers $ 5^{5/4} 6^{4/5}$
is \emph{divisible} by $\delta_L$, and the ratio has
nontrivial $5$-adic valuation.
From the discriminant bounds of Odlyzko \cite{Odlyzko}, under
the assumption of GRH,
one concludes that $[L:\Q] < 2400$ and thus $[L:K] < 24$.
In particular, $L/\Q$ is a solvable extension, and thus
we can apply tools from class field theory. 

\

\Rem  The methods of Odlyzko   \cite{Odlyzko} give 
a general technique to bound the root discriminant of a number
field $K/\Q$. One obtains an estimate of
the form $\delta_K \ge
\gamma(n)$, where $n=[K:\Q]$ and $\gamma(x)$ is
an explicitly calculable continuous convex function
of $x$. Under the assumption of the Generalized
Riemann Hypothesis, one obtains similar bounds, 
with $\gamma$ replaced by some larger (still convex)
function $\gamma'(x)$. Both $\gamma$ and
$\gamma'$ have finite limit as $x$ goes to $\infty$.
In fact, 
$$\lim_{x \longrightarrow \infty} \gamma(x) = 4 \pi e^{\gamma} = 22.38,
\qquad  \lim_{x \longrightarrow \infty} \gamma'(x) \sim 41.$$
In particular, without the GRH, the estimate $\delta_L < 31.349$
does not imply the degree $[L:\Q]$ is bounded, which might
allow us to compute $L$ explicitly.
In the calculations of Fontaine, Schoof and Brumer--Kramer,
all root discriminant estimates fall within the
unconditional (on GRH) bounds of Odlyzko and
so allow the use of unconditional estimates.

\

Our calculations in this section could be shortened by more 
reliance on computer calculation. However, for exposition we include
as much class field theory as we can do by hand. This leads us to
consider several group theory lemmata which allow us to do computations
in smaller fields.

\

The root discriminant of $K$ is
$\delta_K = 5^{23/20} 6^{4/5}$. Since the exponents of
$2$ and $3$ in $\delta_K$ equal those for $\delta_L$,
the extension $L/K$ is unramified outside primes lying 
above $5$.
Let $F = \Q(\zeta_5,576^{1/5})$ ($F$ can also
be written as $\Q(\zeta_5,18^{1/5})$ or
$\Q(\zeta_5,24^{1/5})$). 
 Then $F/\Q(\zeta_5)$ is unramified
at $1 - \zeta_5$. The prime $1 - \zeta_5$ splits completely in $F$.
If $\pi_{F,i}$ for $i=1,\ldots,5$ are the primes above $1-\zeta_5$,
 then
$ N_{F/\Q}(\pi_{F,i}) = 5$.
The extension $K/F$ is totally ramified at all primes $\pi_{F,i}$.
Equivalently,
$\pi_{F,i} = \pi_{K,i}^5$ for all $i$, and $N_{K/\Q}(\pi_{K,i}) = 5$.
Let us consider the factorization of $\pi_{K,i}$ in $L$. Since
$L/\Q$ is Galois, the ramification exponents are equal for all
$i$.  Thus we may write
$$\pi_{K,i} = \prod_{j=1}^{r_{L/K}} \frak{p}^{e_{L/K}}_{i,j},
\qquad N_{L/K}(\frak{p}_{i,j}) = \pi^{f_{L/K}}_{K,i},
\qquad N_{L/\Q}(\frak{p}_{i,j}) = 5^{f_{L/K}},$$
$$r_{L/K} e_{L/K} f_{L/K} = [L:K].$$

\subsection{$L/K$ Tame.}

In this section we assume that $L/K$ is a \emph{tame} extension, and
prove that $L = K$.
\begin{lemma}  \label{lemma:tame} $[L:K] < 10$.
\end{lemma}

\Proof Since $L/K$ is tame,
$\Di_{L/K} =  \prod_{i=1}^5 
\prod_{j=1}^{r_{L/K}} \frak{p}^{e_{L/K} - 1}_{i,j}$,
where $\Di_{L/K}$ is the different. Thus
 $$\Delta_{L/K} = N_{L/K}(\Di_{L/K})
 = \prod_{i=1}^5 \pi^{r_{L/K} f_{L/K}
(e_{L/K}-1)}_{K,i}.$$
Since $N_{K/\Q}(\pi_{K,i}) = 5$,
$\ord_5(N_{K/\Q}(\Delta_{L/K})) = 5 [L:K](1 - 1/e_{L/K})
< 5[L:K]$.
Using the transitivity property of the discriminant
\cite{Serre} we find
$$\delta_{L} = \delta_{K} \cdot N_{K/\Q}(\Delta_{L/K})^{1/[L:\Q]}
< \delta_K \cdot  5^{5/[K:\Q]}
= 5^{23/20} 6^{4/5} 5^{5/100} =  28.925.$$
If $[L:\Q] \ge 1000$, then assuming the GRH,
$\delta_L > 29.094$ by \cite{Odlyzko}. This is a contradiction to
the above upper bound, so 
$[L:K] < 10$. $\qed$

\

\Rem If $5 \  |  \ [L:K]$ and $[L:K] < 10$ then
$[L:K] = 5$ and $\Gal(L/K)$ is tame if and only if it is unramified.
We shall consider the case $L/K$ unramified of degree $5$ in 
section~\ref{sec:unramified}.
Therefore we assume for now that $[L:K]$ has order coprime to $5$.

\begin{lemma} \label{lemma:group1}
Let $\G$ be a finite group, and let $\G' = [\G,\G]$ be its
commutator subgroup. Suppose moreover that $\G/\G' \cong \Z/5 \Z$,
and
$|\G'| \ge 10$. Then $\G \cong \Z/5 \Z$.
\end{lemma}

\Proof Assume that $|\G'| < 10$. It suffices to show that
$\G$ is abelian, since then $\G' = \langle 1 \rangle$ and
$\G \simeq \Z/5\Z$.
If $\G'$ has order $5$, then $\G$ has order $5^2$ and is
necessarily abelian. 
 Thus $\G'$ has order coprime to $5$, and
thus there exists a section $\Z/5\Z \rightarrow \G$. The
action of $\Z/5\Z$ on $\G'$  by conjugation induces an
automorphism of $\G'$.
 Since for all groups $\G'$ of order less than $10$,
$|\Aut(\G')|$ is coprime to $5$, this action must be trivial,
and thus $\G$ is abelian. $\qed$

\begin{lemma} \label{lemma:Htame}
If the extension $L/K$ is  tame and  of degree coprime to $5$,
then $L = K$.
\end{lemma}
\Proof Let $J$ be the field $\Q(\zeta_5,2^{1/5})$.
We have the following exact sequence of groups
$$0 \lra \Gal(L/K) \lra \Gal(L/J) \lra  \Z/5 \Z \lra 0.$$
The extension $L/J$ is Galois since $L/\Q$ is Galois.
By Lemma~\ref{lemma:tame}, $[L:K] < 10$. Thus
by Lemma~\ref{lemma:group1}, either $L = K$ or
$\Gal(L/K)$ is not
the commutator subgroup of $\Gal(L/J)$.
Since by assumption $[L:K]$ has order coprime to $5$,
either $L = K$ (in which case we are done) or
$\Gal(L/J)^{ab}$ is  not a $5$-group. Hence
$J$ admits a
Galois extension $E/J$ contained in
$L$, and of order coprime to $5$. 

\begin{sublemma} $E/J$ is unramified at $2$ and $3$.
\end{sublemma}

\Proof Let $p \in \{2,3\}$. Consider ramification degrees $e_p$.
One has
$$e_p(E/J) \ |  \ e_p(L/J) = e_p(L/K) e_p(K/J).$$
Moreover, $L/K$ is unramified at primes above $2$ and $3$, and
$[K:J] = 5$. Thus $e_p(E/J)$ is either $1$ or $5$, and thus
must be $1$, since $5 \nmid [E:J]$. $\qed$

\

\noindent \bf (Continuation of Proof of Lemma~\ref{lemma:Htame})\rm\ 
 Thus  $E/J$ is a
non-trivial
abelian extension of degree coprime to $5$, unramified outside 
the unique prime $\pi_J$ above $5$, and tamely ramified at
$\pi_J$. 
Such extensions are classified by class field theory.
One has by {\tt pari} that $\Cl(\Of_J) = 1$. On the other hand,
$J/\Q$ is totally ramified at $5$, and so
$(\Of_J/\pi_J \Of_J)^* \simeq \F^*_5$
which  is generated by the global unit
$(1 + \sqrt{5})/2 \equiv -2$. Thus 
$E$ does not exist. This proves that $L = K$. $\qed$

\subsection{$L/K$ Wild.}

Recall that $[L:K] < 24$.
In this section, we assume that $L/K$ is wildly ramified and of
degree $10$, $15$ or $20$, and leave $[L:K] = 5$ until
the next section.
\begin{lemma} \label{lemma:group2}
Let $H$ be a group of order $10$, $15$ or $20$.
Let $\G$ be an extension of $\Z/5 \Z$ by $H$. Then $\G^{ab}$
is not a $5$-group.
\end{lemma}
\Proof  Let $H' \simeq \Z/5\Z$ be a $5$-Sylow subgroup of $H$. Since
$5(1+5) > 20$, $H'$ is the unique $5$-Sylow subgroup of $H$.
Thus $H'$ is preserved under automorphisms of $H$.
Conjugation on $H$ by an element of $\G$ is an automorphism, and thus
$g H' g^{-1} = H'$ for any $g$ in $\G$. In particular,
$H'$ is normal in $\G$.
The quotient group $\G/H'$ has order equal to $|H|$. Thus
$\G/H'$, like $H$, has a normal $5$-Sylow subgroup
and hence has
an abelian quotient of order $2$, $3$ or $4$.
Since this quotient is also a quotient of $\G$,
the lemma is proved. $\qed$ 

\

Let $J$ be the field $\Q(\zeta_5,2^{1/5})$. We have the following
exact sequence of groups
$$0 \lra \Gal(L/K) \lra \Gal(L/J) \lra \Gal(K/J) \lra 0.$$
By  Lemma~\ref{lemma:group2}, $\Gal(L/J)^{ab}$ is not a $5$-group.
Thus $J$ admits an abelian extension of degree coprime to $5$,
contained in $L$.
The non-existence of such an extension was proved during
the proof of  
Lemma~\ref{lemma:Htame}.

\subsection{$L/K$ of degree $5$.}
\label{sec:unramified}

To finish the proof of Theorem~\ref{theorem:5power},
we must consider the cases where $L/K$ is either wildly ramified
of degree $5$, or unramified of degree $5$.
Assume that we are in one of these cases; we shall
arrive at a contradiction.
$\Gal(L/\Q(\zeta_5))$ is a group of order
$125$ that surjects onto $\Z/5 \Z \oplus \Z/5 \Z$. There are three
groups up to isomorphism with this property. All of them admit
at least one morphism
 to $\Z/5 \Z$ with kernel $\Z/5 \Z \oplus \Z/5\Z$ that
factor through the map to $\Gal(K/\Q(\zeta_5))$.
Thus there exists a field $D/\Q(\zeta_5)$, contained within
$K$, such that $\Gal(L/D) \simeq \Z/5 \Z \oplus \Z/5 \Z$.

\begin{lemma} There exists an intermediate field $E$
contained in $L$ and containing $D$ such that $E$ is not 
equal to $K$ and $E/D$ is unramified at primes above $2$ and $3$.
\end{lemma}

\Proof Since the root discriminant of $L$ locally at $2$ and $3$ is bounded by
$2^{4/5}$ and $3^{4/5}$ respectively, this lemma is obvious if the
root discriminant for $D$  attains these bounds, since then any
subgroup of $\Gal(L/D) = \Z/5 \Z \oplus \Z/5 \Z$ not
corresponding to $K$ will produce the required $E$.
Thus we may assume that
$D = \Q(p^{1/5},\zeta_5)$ with $p$ equal to $2$ or $3$. Assume
that $p = 2$. 
The inertia
group $\I_3(L/D)$  is of order $5$, since 
$e_3(L/K) = 1$ and $e_3(K/D) = 5$.
Thus we see that the fixed field $E$ of
$\I_3(L/D) \subset \Gal(L/D)$
  \emph{is} unramified at $3$
above $D$.  Moreover, $E/D$ is unramified above $2$, since $L/K$
and $K/D$ are. Finally, 
$E$ is not $K$ since $K/D$ is ramified at $3$.
An identical argument works for $p = 3$. $\qed$

\begin{lemma} \label{lemma:con} If $D/\Q$ is wildly ramified at $5$ then either
$E/D$ is unramified at $5$ or
$\Delta_{E/D} = \pi^8_D$ where $\pi_D$ is the unique prime above
$5$ in $D$. If $D$ is tamely ramified at $5$, then
 $D = \Q(\zeta_5,24^{1/5})$ and  
$\Delta_{E/D}$ divides $(\pi_{D,1}
\ldots \pi_{D,5})^8$, where $\pi_{D,i}$ are the primes in
$D$ above $5$.
\end{lemma}

\Proof Either $D/\Q$ is tamely ramified at $5$ (in which case it is
$\Q(\zeta_5,24^{1/5})$) or it is not.
Suppose first that $D/\Q$ is wildly ramified.  We
may assume also that $E/D$ is  wildly ramified, since otherwise
it is unramified, and we are done. 
Let $v_{E/D}$ denote the exponent of $\pi_D$ in 
$\Delta_{E/D}$. Suppose that $v_{E/D} \ge 10$. Then
$N_{D/\Q}(\Delta_{E/D}) = 5^{v_{E/D}}  \ge 5^{10}$, so
$$\delta_{E,5} = \delta_{D,5} N_{D/\Q}(\Delta_{E/D})^{1/100}
\ge  5^{23/20} 5^{10/100} = 5^{5/4}.$$
Since $\delta_{E,5}$ divides $\delta_{L,5}$, this estimate
contradicts the bound given by Fontaine
(Theorem~\ref{theorem:font}).
Hence $v_{E/D} < 10$.
On the other hand, 
we have the following equality regarding the different
(\cite{Serre}, IV. Prop.~4):
$$\sum_{i=0}^{\infty} (|\G_i| - 1) = v_{\pi_E}(\Di_{E/D})$$
where $\pi_E$ is the prime in $E$ above $\pi_D$, and $\G_i$ is
the higher ramification group (in the lower numbering) at
$\pi_E$.
Thus if $v_{E/D}$ is the exponent of
$\pi_E$ in the different $\Di_{E/D}$ (equivalently, the exponent
of $\pi_D$ in the discriminant $\Delta_{E/D}$),
 $$v_{E/D} \equiv 0 \mod 4$$
since $|\G_i| - 1$ is either equal to $4$ or $0$.
Thus since $v_{E/D} < 10$, we have
 $v_{E/D} = 4$ or $8$. Since we have wild ramification,
$v_{E/D} >  e_{E/D} - 1$. Thus $v_{E/D} = 8$, and
$\Delta_{E/D} = \pi^8_D$.

Suppose now that
$D= \Q(\zeta_5,24^{1/5})$. Let $\pi_{K,i}$ be the unique prime
above $\pi_{D,i}$ in $\Of_K$. Let 
 $\Delta_{L/K} = (\pi_{K,1}
\ldots \pi_{K,5})^v$ (note all exponents are equal since $L/\Q$
is Galois). 
If $v \ge 10$ then
$$\delta_{L,5} = \delta_{K,5} N_{K/\Q}(\Delta_{L/K})^{1/500}
\ge  5^{23/20} 5^{50/500} = 5^{5/4}$$
contradicting  the bound of Fontaine as above. 
Note that the only ramification in $L/D$ occurs at primes
above $5$. 
Thus  since $v < 10$, 
$$\Delta_{L/D} = \Delta^5_{K/D} N_{K/D} (\Delta_{L/K}) \ | \ 
 (\pi_{D,1} \ldots \pi_{D,5})^{40 + 10}$$
where the left hand side \emph{strictly} divides the
right  as ideals in $\Of_D$.
If $\pi_{D,i}$ occurs in $\Delta_{E/D}$ with exponent $v_{E/D}$ then
$$  \Delta_{L/D} > \Delta^5_{E/D} =  \pi^{5v_{E/D}}_{D,i}$$
and thus $v_{E/D}  < 10$. 
Since $[L:K] = 5$, all ramification groups have
order $5$ or $1$, so as above, $4 | v_{E/D}$ and thus
$v_{E/D} = 0$ or $8$ for each prime $\pi_{D,i}$. In particular
$\Delta_{E/D}$ divides $(\pi_{D,1}
\ldots \pi_{D,5})^8$. 
 $\qed$ 

\

Since $E/\Q$ is not necessarily Galois,
we can not infer in the proof above
that $v_{E/D}$ is the same for each $\pi_{D,i}$.
However, that will not be  necessary for the sequel.

\

Let us recall now the conductor-discriminant
formula (see for example \cite{Neut}, 11.9, p.~557).
\begin{lemma}[Conductor-Discriminant Formula]
 \label{lemma:cdf} Let $B/A$ be an abelian extension
of algebraic number fields.
Then
$$\Delta_{B/A} = \prod_{\chi} \frak{f}(\chi)$$
where the product runs over all characters of
$\Gal(B/A)$.
\end{lemma}

\begin{lemma} \label{lemma:conduct}
If $E/D$ is ramified, and $D/\Q$ is wildly ramified
at $5$ then
the conductor $\frak{f}_{E/D}$ is equal to $\pi^2_D$. If $D/\Q$
is tamely ramified (and so $D = \Q(\zeta_5,24^{1/5})$),
 then the conductor $\frak{f}_{E/D}$ divides
$(\pi_{D,1} \ldots \pi_{D,5})^2$.
\end{lemma}

\Proof The four nontrivial characters
of $\Gal(E/D) \simeq \Z/5 \Z$ are faithful, and 
have conductor $\frak{f}_{E/D}$. Thus 
by the conductor-discriminant formula
(Lemma~\ref{lemma:cdf}), $\Delta_{E/D} = (\frak{f}_{E/D})^4$.
The result then follows from Lemma~\ref{lemma:con}. $\qed$

\

By Lemma~\ref{lemma:conduct}, the possibilities for $E$ 
will be constrained by
the ray class group
of $\frak{f}:=\pi^2_D$ or $(\pi_{D,1} \ldots \pi_{D,5})^2$.
 We may calculate these groups with
the aid of {\tt pari}. The results are tabulated in the table in
the appendix (section~\ref{sec:ray}), and
they indicate the proof of
Theorem~\ref{theorem:5power}
is complete,
after noting that all cases where the ray class field is
nontrivial, $K$ is the ray class field, contradicting the
definition of~$E$.

\section{$N = 10$.}

Let us begin by stating the analogs for $N=10$
of Theorems~\ref{theorem:5power} and \ref{theorem:Schoof}
(another special case of Theorem
3.3 and the proof of Corollary 3.4 in \cite{Schoof}):

\begin{theorem} \label{theorem:3power} 
Let $G/\Z[\frac{1}{10}]$ be a finite group scheme of
$3$-power order such that:
\begin{enumerate}
\item  Inertia at $2$ and $5$ acts through
a cyclic $3$-group.
\item  The action of Inertia on the subquotients
$G[3^n](\Qbar)/G[3^{n-1}](\Qbar)$ is through an order $3$ quotient for all $n$.
\end{enumerate}
 Assume the GRH discriminant bounds of
Odlyzko.
Then $G$ has a filtration by
the group schemes $\Z/3 \Z$ and $\mu_3$. Moreover, if $G$ is
killed by $3$, then $\Q(G) \subseteq H$, where
$K:= \Q(\sqrt[3]{2},\sqrt[3]{5},\zeta_3)$, and $H$ is the
Hilbert class field of $K$, which is of degree $3$ over $K$.
\end{theorem}

\begin{theorem}[Schoof] \label{theorem:schoof2}
Let $p = 2$ or $5$.
Let $G/\Z[\frac{1}{p}]$ be a finite group scheme of
$3$-power order such that inertia at $p$ acts through
a cyclic $3$-group. Then $G$ has a filtration by
the group schemes $\Z/3 \Z$ and $\mu_3$. Moreover, the extension
group $\Ext^1(\mu_3,\Z/3 \Z)$ of group schemes over $\Z[\frac{1}{p}]$ is
trivial, and there exists an exact sequence of group schemes
$$0 \lra M \lra G \lra C \lra 0$$
where $M$ is a diagonalizable group scheme over $\Z[\frac{1}{p}]$,
 and $C$ is
a constant group scheme.
\end{theorem}

We delay the proof of theorem~\ref{theorem:3power} until
Section~\ref{sec:porf}.

\

The situation for $N=10$ is analogous to $N=6$, but 
one technical difficulty is that 
$\I_p(H/\Q)$ is not a normal subgroup of $\Gal(H/\Q)$, for
either $p = 2$ or $p = 5$.
We do however prove the following:
\begin{lemma}  \label{lemma:splitting}
The primes $2$ and $5$
split into exactly $3$ distinct primes in $H$.
\end{lemma}
\Proof  The equivalent statement is true with $H$ replaced by
$K$. Thus it suffices to prove that $2$ and $5$ remain inert
in the extension $H/K$.
The easiest way to see this is by 
comparing residue field degrees.
$H$ is the compositum of $K$ and the Hilbert
class field $B$ of $\Q(20^{1/3})$. In  $\Q(20^{1/3})$, the primes above
$2$ and $5$
are not principal (a simple check with {\tt pari}), and so remain
inert in $B$.
 Thus for $p =2$ and $p=5$ the residue field degree
$f_{B/\Q} = 9$. Since for $p=2$ and $p =5$ one
finds that $f_{K/\Q} = 3$, it follows that the
the primes above $2$ and $5$ in $K$
remain inert in $H$. $\qeda$ 

\

It follows from Lemma~\ref{lemma:splitting} that
the subgroups $\D_p(H/\Q)$ are of index three in $\Gal(H/\Q)$.
One sees that $\I_{p'}(H/\Q)$ is a set of representatives
for the cosets of $\D_p(H/\Q)$ in $\Gal(H/\Q)$, where
$\{p,p'\} = \{2,5\}$ (as an unordered pair), since the
possible fixed fields of $\D_p(H/\Q)$ are never fixed
by $\I_{p'}(H/\Q)$.
This leads to the following construction:

\begin{lemma} \label{lemma:once}
 Let  $\{p,p'\} = \{2,5\}$. Let $\ovM \subset A[3]$ be
a $\D_p(H/\Q)$-module. Let $\sigma \in 
\I_{p'}(H/\Q)$ be a nontrivial element.
Then
$$\haM = \ovM + \sigma \ovM = \ovM + (\sigma-1) \ovM$$
is a $\Gal(\Qbar/\Q)$-module. Moreover, 
$\dim(\haM) \le 2 \dim(\ovM)$  with equality
if and only if  $\ovM$ and
$(\sigma - 1) \ovM$ have trivial intersection inside $A[3]$.
\end{lemma}

\Proof By Grothendieck (Theorem~\ref{theorem:Groth})
 one finds that as an endomorphism, for $\sigma \in \I_{p'}$,
we have $(\sigma -1)^2 = 0$ on $A[3]$. Thus $\haM$ is an $\I_{p'}$-module.
On the other hand, $\ovM$ is a $\D_p$-module, and
the coset space $\Gal(H/\Q)/\D_p(H/\Q)$ is given by
$\I_{p'}(H/\Q)$. Thus $\haM$ is a $\Gal(H/\Q)$-module.
The final claim is clear from the construction of
$\haM$.
$\qed$.

\

We now apply this construction not to $\ovMt(p)$, as in 
section~\ref{sec:const},
but to $\ovMf(p)$.

\begin{lemma}  \label{lemma:extrareview}  Fix 
$p \in \{2,5\}$ and assume that $\ord_3(\Phi_{\haA}(p))$ is
maximal. Then 
$$A[3] = \haMf(p) + \ovMt(p).$$
\end{lemma}

 If $\kappa = \haMf(p)$, then from Theorem~\ref{theorem:BK},
$$\ord_{3}(\Phi_{\haA'}(p)) - \ord_{5}(\Phi_{\haA}(p))
=  \dim \ \kappa \cap \ovMf(p) + \dim \ \kappa \cap \ovMt(p)
- \dim \ \kappa.$$
Since $\ovMf(p) \subseteq \kappa \cap \ovMt(p)$,
we find that this quantity is at least $2 t_p - \dim \ \kappa$. 
By the maximality assumption on $\ord_3(\Phi_{\haA}(p))$ we conclude
that $2 t_p - \dim \ \kappa \le 0$.
On the
other hand, by Lemma~\ref{lemma:once}
 we see that $\dim \kappa \le 2 t_p$,
with equality if and only if $\ovM$ and $(\sigma-1) \ovM$ have 
trivial intersection inside $A[3]$. Thus $\dim \kappa = 2 t_p$,
and $\ovM$ and $(\sigma -1) \ovM$ have trivial intersection
inside $A[3]$.
By Lemma~\ref{lemma:sigmaminus},
the image of $(\sigma - 1)$ on $A[3]$
for $\sigma \in \I_{p'}$  is contained within $\ovM_{2}(p')$, and
thus
has dimension at most $t_{p'}$. This immediately proves that
$t_p \le t_{p'}$, and by symmetry, that $t_2 = t_5$.
Moreover, equality also forces $\ovMf(p) = \kappa \cap \ovMt(p)$.
In particular, $\dim(\ovMt(p) +  \haMf(p)) = 
\dim(\ovMt(p)  + \kappa)$ is at least
$$\dim \ovMt(p)  + \dim \kappa - \dim \kappa \cap \ovMt(p)
= \dim \ovMt(p)  + \dim \kappa - \dim \ovMf(p)$$
which equals $(2d - t_p) + 2 t_p - t_p = 2d$.
In other words,
$$A[3] = \haMf(p) + \ovMt(p).$$
$\qed$.

\begin{lemma}  Let  $\{p,p'\} = \{2,5\}$. \label{lemma:unram}
 For $\ord_3(\Phi_{\haA}(p))$ maximal, $\Q(A[3])$ is unramified
at $p$.  \end{lemma}
\Proof From Lemma~\ref{lemma:extrareview},
we have  that  $A[3] = \haMf(p) + \ovMt(p)$. 
Moreover, we have a direct sum decomposition
$$\haMf(p) = \ovMf(p) \oplus (\sigma-1) \ovMf(p)$$ 
which follows from the discussion above, since for
$\haMf(p)$ to have dimension $2t_p$, 
$\ovMf(p)$ and $ (\sigma-1) \ovMf(p)$
must have trivial intersection.
By definition,
$\I_p$ acts trivially on $\ovMt(p)$. Thus it suffices to show that
$\I_p$ acts trivially on $\haMf(p)$.
Since $\I_p(H/\Q) = \I_p(H/\Q(\zeta_3))$ for $p \in \{2,5\}$ it
clearly suffices to show that
$\I_p(H/\Q(\zeta_3))$ acts trivially on $\haMf(p)$ considered
as a $\Gal(H/\Q(\zeta_3))$-module.

\

Fix a basis for $\ovMf(p)$, and choose the corresponding
basis for $(\sigma-1) \ovMf(p)$ whose basis elements
are $(\sigma-1)$ of our chosen basis for $\ovMf(p)$.
Note that $(\sigma-1)$ \emph{fixes} $(\sigma-1) \ovMf(p)$
since $(\sigma-1)^2 = 0$. Thus with 
respect to this basis, the action of $\sigma$ on
$\haMf(p)$ is given by
the following matrix
$$\rho(\sigma)
=\left( \begin{array}{cc} \Id_t & 0 \\ \Id_t & \Id_t \end{array}
\right)$$
where $\Id_t$ is an element of $M_t(\F_p)$, the $t \times t$
matrices over $\F_p$, and $\rho(\sigma)$ denotes the image of
$\sigma$ in $\Gal(\Q(\haMf(p))/\Q(\zeta_3))$. 
By Lemma~\ref{lemma:sigmaminus},
for $\tau \in \I_p$,  the module  $(\tau - 1)A[3]$ lies within
$\ovMf(p)$. Moreover, $\ovMf(p) \subseteq \ovMt(p)$
is fixed by $\tau$.  Thus with respect to our basis
the action
of  $\tau \in \I_p$ 
is given by a matrix
$$\rho(\tau)
 =\left( \begin{array}{cc} \Id_t & a \\ 0 & \Id_t \end{array} \right)$$
for some $a \in M_t(\F_p)$.
It suffices to prove that $a = 0$, since then we have shown
$\I_p$ acts trivially
on $A[3]$.
Since 
$$\Q(\haMf(p)) \subseteq
\Q(A[3]) \subseteq H$$
and since $[H:\Q(\zeta_3)] = 27$,  the
group $U$ generated by $\langle \rho(\sigma), \rho(\tau) \rangle$ 
is a group of order
dividing $27$. Moreover, since all powers of $a$ commute
with one another,
$U$ is a subgroup
of $\GL_2(\F_3[a])$.
Thus our desired conclusion
($a = 0$)
follows from the following result:

\begin{sublemma}  Let $U \subseteq \GL_2(\F_3[a])$ be a subgroup
of order dividing $27$ containing the elements
$$\left< \left( \begin{array}{cc} 1 & a \\ 0 & 1 \end{array} \right),
\left( \begin{array}{cc} 1 & 0 \\ 1 & 1 \end{array} \right) \right>$$
then $a = 0$.
\end{sublemma}

\Proof If $U$ is abelian, then the
commutator $[\sigma,\tau] = 1$, and one computes
immediately that $a = 0$. 
 Thus $U$ has order $27$. From a classification of all non-abelian
groups of order $27$, we find that $[U,U]$ is  central.
From
$[\sigma,\tau] \sigma - \sigma [\sigma,\tau] = 0$, we infer
that
that $a^2 = 2a + a^2 = 0$. In characteristic $3$, this proves that
$a = 0$. $\qed$

\

\Rem The utility in Lemma~\ref{lemma:unram} is that it allows
us to reduce arguments for $N = 10$ to the exact
analogs of  arguments
for $N =6$. The reason is as follows.
When  $\ord_{3}(\Phi_{\haA}(p))$
is maximal, Lemma~\ref{lemma:unram} implies that $\Q(A[3])$
is unramified at $p$. Thus $\Q(A[3])$ must be a Galois
(over $\Q$) 
subfield
of $H$ fixed by all conjugates  of the inertia
group $\I_p(H/\Q)$.
The only such fields are the
subfields of $\Q(\zeta_3,p'^{1/3})$.  
Thus we  obtain the desired results of
section~\ref{sec:const}
\emph{without} the hypothesis that $\I_p(H/\Q)$
is normal in $\Gal(H/\Q)$.

\

We may now establish Theorem~\ref{theorem:results5}
in much the same way as Theorem~\ref{theorem:results}. 
 Let  $\{p,p'\} = \{2,5\}$.
Here are the extra steps required 
to complete the proof:
\begin{enumerate}
\item For $\ord_{3}(\Phi_{\haA}(p))$ maximal, $\Q(A[3])$
is unramified at $p$ by Lemma~\ref{lemma:unram}.
Thus $A[3]$  prolongs to a finite group
scheme over $\Z[1/p']$, and
from Theorem~\ref{theorem:schoof2} we infer there exists
an exact sequence of group schemes (over $\Z[1/p']$)
$$0 \lra  \mu^m_3 \lra A[3] \lra (\Z/3 \Z)^n \lra 0.$$
The arguments of Lemma~\ref{lemma:newreview}  apply
\emph{mutatis mutandis} to show that $m=n=d$ and $A$ has ordinary reduction at $3$.
\item The arguments of section~\ref{sec:nontoric} 
and Lemma~\ref{lemma:mult} hold with $5$
replaced by $3$ and ``$2$ and $3$'' replaced by
``$2$ and $5$''. It suffices
to show that $\Q$ does not admit any Galois extensions
of order $3$ unramified outside $2$ and $5$. This
follows from the Kronecker--Weber theorem.
\item The maximal Galois extension of $\Q$,
contained within  $H$  and unramified at
$2$ is $\Q(\zeta_3,5^{1/3})$. Hence the proof of
Lemma~\ref{lemma:hat3} still applies, since the inertia
subgroups of $\Gal(\Q(A[3])/\Q)$ when
$\Q(A[3]) \subseteq \Q(\zeta_3,5^{1/3})$ are \emph{normal}.
Similarly, a proof of Lemma~\ref{lemma:above} requires us only to
note that
the maximal 
unramified at $3$ extension of $\Q(\zeta_3)$ inside $H$ is
$\Q(\zeta_3,10^{1/3})$, which \emph{is} ramified at $5$.
\item A final contradiction is reached because
$$3^{4d} \le (1 +\sqrt{3})^{4d}$$ 
is not true for $d > 0$. One might remark at this point that since
$A$ has good reduction
at $3$, and since $A$ is defined over $\Q$, the 
rational $3$-torsion
injects into $A(\F_p)[3^{\infty}]$, as follows from standard facts
about formal groups (see, for example,
\cite{KatzII}).  Note this fact is not essential,
however, since for any prime $q$ of good reduction,
an inclusion of a constant group scheme
$C$ into $A$ \emph{automatically} reduces to an
inclusion $C(\F_q) \hookrightarrow A(\F_q)$.
\end{enumerate}
Thus it remains to prove Theorem~\ref{theorem:3power}.

\subsection{Group Schemes over $\Z[1/10]$.}

\label{sec:porf}  
Since $\Gal(H/\Q(\zeta_3))$ is a $3$-group, the discussion at the
beginning of section~\ref{sec:groups} shows that it suffices to
prove that if
$L = \Q(G \times G_{-1} \times G_{2} \times G_{5})$ then
$L \subseteq H$. 
One has the following estimate
of the root discriminant for $L$
$$\delta_L <  3^{1 + \frac{1}{3-1}} 2^{1 - \frac{1}{3}}
5^{1 - \frac{1}{3}} = 3^{3/2} 10^{2/3} = 24.118  < 24.258.$$
From the estimates of \cite{Odlyzko} one finds that
$[L:\Q] < 280$, and so $[L:K] < 16$.
One sees that
that $K:=\Q(\sqrt{-3},\sqrt[3]{2},\sqrt[3]{5}) \subseteq L$.
We wish to prove that $\Gal(L/K)$ is a $3$-group.
The root discriminant of $K$ is
$\delta_K = 3^{7/6} 10^{2/3}$,
 and so $L/K$ is at most ramified at primes above $3$.
Let $F = \Q(\sqrt{-3},\sqrt[3]{10})$. Then
$F/\Q(\sqrt{-3})$ is unramified at $\sqrt{-3}$. The prime
$\sqrt{-3}$ splits completely in $F$, and we write
$(\sqrt{-3}) = \pi_{F,1} \pi_{F,2} \pi_{F,3}$,
 $N_{F/\Q}(
\pi_{F,i}) = 3$.
The extension $K/F$ is totally ramified at each $\pi_{F,i}$. One
has
$\pi_{F,i} = \pi_{K,i}^3$ for all $i$, and $N_{K/\Q}(\pi_{K,i}) = 3$.
\subsection{$L/K$ Tame}

In this section we assume that $L/K$ is a \emph{tame} extension.
\begin{lemma}   One has $[L:K] \le 6$.
\end{lemma}
\Proof Arguing as in Lemma~\ref{lemma:tame}
we find that $N_{K/\Q}(\Delta_{L/K}) < 3^{3 [L:K]}$. Thus
$$\delta_L = \delta_K N_{K/\Q}(\Delta_{L/K})^{1/[L:\Q]} \le
3^{7/6} 10^{2/3} 3^{3/18} = 20.082.$$
Yet from the GRH Odlyzko bound, if $[L:\Q] \ge 126$, then
$\delta_L > 20.221 $. Thus we find $[L:K] \le 6$. $\qed$

\

If $[L:K] \le 6$, then either $\Gal(L/K)$ is a $3$-group or
it surjects onto a nontrivial group of order coprime to $3$.
In the latter case, $L$ would contain an abelian extension $J/K$ tamely
ramified and of degree coprime to $3$. 
Moreover, since $L/K$ is unramified at
primes above $2$ and $5$, a similar conclusion holds for $J$.

\begin{lemma} There are no abelian extensions $J/K$ tamely ramified 
at primes above $3$ in $K$, of
order coprime to $3$, and unramified outside $3$. 
\end{lemma}

\Proof We proceed via class field theory. According
to {\tt pari}, the class number of
$K$ is $3$, its Hilbert class field $H$ being the compositum of
$K$ and the Hilbert class field of $\Q(\sqrt{-3},\sqrt[3]{20})$.
Thus by class field theory 
it suffices  (since $J/K$ has order prime to $3$)
to show that global units of $\Of_K$
 generate $(\Of_K/\pi_{K,1} \pi_{K,2} \pi_{K,3})^*$.
 On the other hand, since $K/F$ is totally
ramified, we have an isomorphism
$$(\Of_K/\pi_{K,1} \pi_{K,2} \pi_{K,3})^*
\simeq (\Of_F/\pi_{F,1} \pi_{F,2} \pi_{F,3})^*
\simeq \F^*_3 \times \F^*_3 \times \F^*_3$$
Hence it suffices to use global units from $\Of_F$.
Let $v = (\sqrt[3]{10} - 1)/\sqrt{-3}$. Then from {\tt pari}
we find that the
$2$ fundamental units of $\Of_F$ are given by
$$\epsilon_1 = \frac{1}{4} v^4 - \frac{1}{2} v^2
+ \frac{3}{2} v - \frac{1}{4}$$
$$\epsilon_2 =  \frac{1}{4} v^4 - \frac{1}{2} v^3
+ \frac{3}{2} v^2  - \frac{1}{4}$$
We find that the images of $-1$, $\epsilon_1$, and $\epsilon_2$
in $\Of_F/\pi_1 \times \Of_F/\pi_2 \times \Of_F/\pi_3$ are
$(-1,-1,-1)$, $(1,1,-1)$ and $(1,-1,1)$ respectively. Since
these elements generate the group  $(\F^*_3)^3$, we are done. $\qed$

\subsection{$L/K$ Wild}

We assume that $L/K$ is wildly ramified at $3$, and
(for the moment) 
not a $3$-group. If $\Gal(L/K)^{ab}$ is not a $3$-group,
then there would exist a corresponding extension $J/K$ tame
of order coprime to $3$. Since no such extensions exist
(see the tame case), we may
also assume that
$\Gal(L/K)^{ab}$ is a $3$-group. Let $\G$ denote the
group $\Gal(L/K)$. There should be no confusion between the
group $\G$ and the group scheme $\G$, which will not appear again.
Let $n = |\G|$. Since $n < 16$, we have
$n \in \{6,12,15\}$. All groups of order $15$ are abelian.
If $n = 6$, the only non-abelian group is $S_3$. Yet
$S^{ab}_3 = \Z/2 \Z$. Thus $n =12$. The only group $\G$ of
order $12$ such that $\G^{ab} = \Z/3 \Z$ is the
alternating group $A_4$, which is a
nontrivial extension of
$\Z/3 \Z$ by $\Z/2 \Z \times \Z/2 \Z$. 
\begin{lemma} One has
$N_{K/\Q}(\Delta_{L/K}) \ge 3^{66}$, and $N_{K/\Q}(\Delta_{L/K})
\le 3^{69}$.
\end{lemma}
Assume the first bound is violated.
Then since $N_{K/\Q}(\Delta_{L/K}) = 3^{3k}$ for some $k$, it must be
bounded by $3^{63}$.
Since $[L:\Q] = 12 \times 18 = 216$,
$$\delta_L = \delta_K N_{K/\Q}(\Delta_{L/K})^{1/[L:\Q]} \le
3^{7/6} 10^{2/3} 3^{63/216} = 23.039 < 23.089.$$
Yet from the GRH Odlyzko bound,
$\delta_L > 23.089 $. The other inequality is violated if and only
if $N_{L/K}(\Delta_{L/K}) \ge 3^{72}$. Yet in this case,
$$\delta_L = \delta_K N_{K/\Q}(\Delta_{L/K})^{1/[L:\Q]} \ge 
3^{7/6} 10^{2/3} 3^{72/216} = 3^{3/2} 10^{2/3}.$$
Since $\delta_L <  3^{3/2} 10^{2/3}$ by the Fontaine bound, this
is a contradiction and thus we are done.
$\qed$ 

\

Denote the primes in $L$ above $\pi_{K,i}$
by $\frak{p}_{i,j}$.
Fix a prime $\frak{p} = \frak{p}_{i,j}$,
 and let $\G_0 \subseteq \G$ (respectively, $\G_1 \subseteq \G$)
be the inertia (respectively wild inertia) subgroup corresponding
to $\frak{p}$. Since $L/\Q$ is Galois, 
these groups are well defined up to conjugation.
Moreover, $\G_1$ is normal in $\G_0$, and $\G_0/\G_1$ has prime
to $3$ order. Since $L/\Q$ is Galois, 
$\G_i$ is determined by $\frak{p}$ up to conjugation.
Let us simplify some notation.
Let $v = v_{\frak{p}}(\Di_{L/K})$, $f = f_{L/K}$,
$e = e_{L/K}$, $r = r_{L/K}$. 
 We have standard equalities
$$\Di_{L/K} = \prod_{i=1}^3 \prod_{j=1}^{r} \frak{p}^{v}_{i,j} \qquad
\Delta_{L/K} = \prod_{i=1}^3 \pi^{f r v}_{K,i}
\qquad N_{K/\Q}(\Delta_{L/K}) = 3^{3 f r v}.$$
From the previous lemma,
$22 \le f r v \le 23$. Moreover, $f r e =[L:K] =12$.
Since $K/L$ is wildly ramified, $3$ divides $e$. Hence it suffices to
show that $e = 3$, $e = 6$  and $e = 12$ all
lead to contradictions. If $e = 3$, then $f r = 4$.
Yet $f r$ divides $23$ or $22$, which is impossible. Suppose that
$e = 6$. Then $|\G_0| = 6$, and $\G_0$  must be a normal subgroup of
$\G$ since it is a subgroup of index $2$.
 If $\G$ had such a subgroup, then $\G^{ab}$ would not be
a $3$-group.  Thus we may assume that $e = 12$. If $e =12$
then the $3$-group $\G_1$ would be a normal subgroup of
$\G_0 = \G$. Since $\G$ has no such subgroup, we are done, and
$\Gal(L/K)$ is a $3$-group.

\

We may therefore assume that
$L/K$ is Galois  of degree dividing $9$, and thus abelian. 
Let  $\frak{f}_{L/K}$ be the conductor of this extension.

\begin{lemma} \label{lemma:pickle}
The conductor $\frak{f}_{L/K}$ divides
$(\pi_{K,1} \pi_{K,2} \pi_{K,3})^2$.
\end{lemma}

\Proof Assume otherwise. Since $L/\Q$ is Galois we infer
that $(\pi_{K,1} \pi_{K,2} \pi_{K,3})^3 | \frak{f}_{L/K}$.
We consider separately the three possible Galois groups.

\

Suppose that $\Gal(L/K) \simeq \Z/9\Z$. Then since
$\Z/9 \Z$ has six faithful characters, by the
 conductor-discriminant
 formula (Lemma~\ref{lemma:cdf}), 
$(\frak{f}_{L/K})^6 | \Delta_{L/K}$.
In particular, if 
 $(\pi_{K,1} \pi_{K,2} \pi_{K,3})^3 | \frak{f}_{L/K}$, then
$(\pi_{K,1} \pi_{K,2} \pi_{K,3})^{18} | \Delta_{L/K}$ and so
$$\delta_{L,3} = \delta_{K,3} N_{K/\Q}(\Delta_{L/K})^{1/[L:\Q]} \ge
3^{7/6} 3^{54/162} = 3^{3/2}$$
contradicting the Fontaine bound.

\

Suppose that $\Gal(L/K) \simeq \Z/3\Z \times \Z/3\Z$.
Let $\chi_i$ for $i = 1,\ldots,4$ be the four characters
of order $3$ corresponding to
the four degree $3$ extensions $K_i$ of $K$  inside $L$.
Since the compositum $K_i.K_j = L$ for $i \ne j$, the
lowest common divisor $(\frak{f}(\chi_i),\frak{f}(\chi_j))$ must equal
$\frak{f}_{L/K}$ for all $i \ne j$. In particular,
$(\frak{f}_{L/K})^2 | \prod_{i=1}^{4} \frak{f}(\chi_i)$, and
thus since there are four non-trivial characters of
$\Z/3\Z \times \Z/3\Z$, by the conductor-discriminant
formula (Lemma~\ref{lemma:cdf}), $(\frak{f}_{L/K})^6 | \Delta_{L/K}$,
and (as above) this contradicts the Fontaine bound.

\

Suppose that $\Gal(L/K) \simeq \Z/3\Z$. Then by the
conductor-discriminant formula (Lemma~\ref{lemma:cdf}),
$(\frak{f}_{L/K})^2 | \Delta_{L/K}$, and so
$$\delta_{L,3} = \delta_{K,3} N_{K/\Q}(\Delta_{L/K})^{1/[L:\Q]} \ge
3^{7/6} 3^{18/54} = 3^{3/2}$$
contradicting the Fontaine bound. $\qeda$

\

By Lemma~\ref{lemma:pickle}, since $L/K$ has conductor
$\frak{f}_{L/K}$ dividing
$\frak{f}:=(\pi_{K,1} \pi_{K,2} \pi_{K,3})^2$, the possible
$L$ are constrained by
the ray class group
of $\frak{f}$. The result is tabulated in the table in
the appendix (section~\ref{sec:ray}).
In particular, we find that the ray class field of
conductor $\frak{f}$ is exactly equal to $H$, the Hilbert class
field of $K$, completing our proof of
Theorem~\ref{theorem:3power}.

\section{Appendix.}

\subsection{Ray Class Fields.}
\label{sec:ray}

Here are some class field computations done
using {\tt pari}. They took between fifteen minutes and an hour each
on a Sparc--II ULTRA machine running at $333$ MHz.
The essentials of the {\tt pari} script are 
below.

\begin{center}
\begin{tabular}{|c|c|c|c|}
\hline
$K$ & $\delta_K$ & $\frak{f}$ & $|\Cl_{\frak{f}}|$\\
\hline
$\Q(\zeta_5,2^{1/5})$ & $5^{23/20} 2^{4/5}$ &
$\pi^2_K$ & $1$ \\
$\Q(\zeta_5,3^{1/5})$ & $5^{23/20} 3^{4/5}$ &
$\pi^2_K$ & $1$ \\
$\Q(\zeta_5,6^{1/5})$ & $5^{23/20} 6^{4/5}$ &
$\pi^2_K$ & $5$ \\
$\Q(\zeta_5,12^{1/5})$ & $5^{23/20} 6^{4/5}$ &
$\pi^2_K$ & $5$ \\
$\Q(\zeta_5,24^{1/5})$ & $5^{3/4} 6^{4/5}$ &
$(\pi_{K,1} \ldots \pi_{K,5})^2$ & $5$ \\
$\Q(\zeta_5,48^{1/5})$ & $5^{23/20} 6^{4/5}$ &
$\pi^2_K$ & $5$ \\
\hline
\hline
$\Q(\zeta_3,2^{1/3},5^{1/3})$ & $3^{7/6} 10^{2/3}$ & $(\pi_{K,1} \pi_{K,2}
\pi_{K,3})^2$ & $3$ \\
\hline
\end{tabular}
\end{center}

\

\subsection{Pari Script.}
Here is the {\tt pari} script for fields other than
$\Q(\zeta_5,24^{1/5})$ and $\Q(\zeta_3,2^{1/3},5^{1/3})$, where
an  adjustment
must be made since the conductor is of a slightly different form. 
The calculation of the discriminant was included as a check against
typographical errors in the defining 
polynomials\footnote{Bjorn Poonen points out to me that as written,
this {\tt pari} program assumes the GRH during computation.
One may remove this assumption from the calculation using
{\tt bnfcertify}.}.

\begin{verbatim}
allocatemem()
allocatemem()
allocatemem()
allocatemem()
nf=nfinit(poly defining K);
factor(nf[3])
bnf=bnfinit(nf[1],1);
pd=idealprimedec(nf,5);
pd1=idealhnf(nf,pd[1]);
idealnorm(nf,pd1)
pd2=idealmul(nf,pd1,pd1);
bnrclassno(bnf,pd1)
bnrclassno(bnf,pd2)
\end{verbatim}

\noindent \it Email address\rm:\tt \  fcale@math.berkeley.edu
\end{document}